\documentclass[12pt]{amsart}
\usepackage{amsmath,amssymb, xypic,verbatim,amscd}
\numberwithin{equation}{section}

\newcommand\bull{\sssize{\bullet}}
\newcommand\home{\operatorname{Hom}}

\newcommand\mf{\mathcal{F}}
\newcommand\mk{\mathcal{K}}

\newcommand\me{\mathcal{E}}
\newcommand\aA{{\Theta} }
\newcommand\mfz{\mathfrak{Z}}
\newcommand\mfa{\mathfrak{A}}

\newcommand\Fl{\operatorname{Fl}}

\newcommand\tensor{\otimes}
\newcommand\ml{\mathcal{L}}

\newcommand\codim{\operatorname{codim}}

\newcommand\GL{\operatorname{GL}}
\newcommand\mj{\mathcal{J}}
\newcommand\Fone{J}
\newcommand\Ftwo{\mj}
\newcommand\mg{\mathcal{G}}
\newcommand\mh{\mathcal{H}}

\newcommand\rk{\operatorname{rk}}

\newcommand\mi{\mathcal{I}}

\newcommand\Gr{\operatorname{Gr}}
\newcommand\SL{\operatorname{SL}}

\newcommand\mm{\mathcal{M}}
\newcommand\muu{\mathfrak{U}}
\newcommand\mukv{\mathfrak{U}_{\mk}(V)}
\newcommand\muiw{\mathfrak{U}_{\mi}(W)}
\newcommand\mhh{\mathfrak{H}}
\newcommand{\leto}[1]{\stackrel{#1}{\to}}
\newcommand\hoo{\home_{\mi}(V,Q,\mf,\mg)}
\newcommand\hoof{\home_{\mi}(V,Q,\me,\mh)}

\newcommand\ham{\mhh_{\mi}(V,Q,\mk)}

\newtheorem{theorem}{Theorem}[section]

\newtheorem{proposition}[theorem]{Proposition}

\newtheorem{lemma}[theorem]{Lemma}

\newtheorem{definition}[theorem]{Definition}
\newtheorem{defi}[theorem]{Definition}

\newtheorem{remark}[theorem]{{ Remark}}

\theoremstyle{remark}

\begin{document}
\title[Fulton's conjecture]{Geometric Proof of a Conjecture of Fulton}
\author{Prakash Belkale}
\address{Department of Mathematics\\ UNC-Chapel Hill\\ CB \#3250, Phillips Hall
\\ Chapel Hill,  NC 27599}
\email{belkale@email.unc.edu}
\begin{abstract}
We give a geometric proof of a conjecture of Fulton on the
multiplicities of irreducible representations in a tensor product of
irreducible representations for $\GL(r)$. This conjecture was proven
earlier by Knutson, Tao and Woodward using the Honeycomb theory.
\end{abstract}
\maketitle

\section{Introduction}
Recall that irreducible {\em polynomial} representations of $\GL(r)$
are indexed by sequences $\lambda=(\lambda_1\geq\dots\geq
\lambda_r\geq 0)\in \Bbb{Z}^r$. Denote the representation
corresponding to $\lambda$ by $V_{\lambda}$. Define
Littlewood-Richardson  coefficients $c^{\lambda}_{\mu,\nu}$ by:
$V_{\mu}\otimes V_{\nu}=\sum c^{\lambda}_{\mu,\nu}V_{\lambda}.$ W.
Fulton conjectured that for any positive integer $N$,
$$c^{\lambda}_{\mu,\nu}=1\Leftrightarrow c^{N\lambda}_{N\mu,N\nu}=1.$$
This conjecture was proved by A. Knutson, T. Tao and C. Woodward
~\cite{KTW} using the Honeycomb theory.

In this article we give a geometric proof of Fulton's conjecture
based the geometric proof of Horn and saturation conjectures given
in ~\cite{gh}. The techniques in the proof to be given  can be
applied in quantum cohomology (this is our main motivation, see
~\cite{qh} for the multiplicative generalization of Horn and
saturation conjectures), and hopefully also in quiver theory, to
prove analogues of Fulton's conjecture.

Our proof deduces Fulton's conjecture from the projectivity of some
Geometric invariant theory (GIT) moduli spaces, a technique which is
sufficiently categorical  for generalizations. This technique is
most easily understood in the geometric proof of Fulton's original
conjecture given here.

I thank Harm Derksen for useful discussions.
\subsection{Conventions}
 We make the following conventions:
\begin{itemize}
\item An integer $s\geq 1$ will be fixed for the proof.
\item For a vector space $W$, let $\Fl(W)$ denote the variety of complete flags on it. If $\me\in \Fl(W)^s$, we will assume that $\me$ is
written in the form $(E^1_{\bull},\dots,E^s_{\bull})$.
\item We use the notation $[n]=\{1,\dots,n\}$.
\end{itemize}

\section{Some results and notation from ~\cite{gh}}
In this section we recall some results and notation from ~\cite{gh}.
The reader may wish to turn to Section ~\ref{shostakovich} now.
\subsection{Schubert cells in Grassmannians}
Let $I\subseteq [n]=\{1,\dots,n\}$ be a subset of cardinality $r$.
Such a set will always be written as $I=\{i_1<\dots<i_r\}$. Let
$$E_{\sssize{\bullet}}:\text{ }\{0\}=E_0\subsetneq E_1\subsetneq\dots\subsetneq E_n=W$$
be a complete flag in an $n$-dimensional vector space $W$. Define
the Schubert cell $\Omega^o_{I}(E_{\sssize{\bullet}})\subseteq
\Gr(r,W)$ by
$$\Omega^o_I(E_{\sssize{\bullet}})=\{V\in \Gr(r,W)\mid \rk(V\cap E_u)=a \text{ for } i_a\leq u< i_{a+1},\text{ }a=0,\dots,r\}$$
where $i_0$ is defined to be $0$ and $i_{r+1}=n$.
$\Omega^o_{I}(E_{\sssize{\bullet}})$ is smooth. Its closure will  be
denoted by  ${\Omega}_{I}(E_{\sssize{\bullet}})$. For a fixed
complete flag on $W$, it is easy to see that (~\cite{ful2}, \S 1)
every $r$-dimensional vector subspace  belongs to a unique Schubert
cell.
\begin{defi}Let $V$ be a $r$-dimensional subspace of an
$n$-dimensional vector space $W$, and $\me\in \Fl(W)^s$. Let
$I^1,\dots, I^s$ be the unique subsets of $[n]$ each of cardinality
$r$ such that $V\in \Omega^o_{I^j}(E^j_{\bull})$ for $j=1,\dots,s$.
Define $\dim(V,W,\me)$ to be the expected dimension of the
intersection $\cap_{j=1}^s\Omega^o_{I^j}(E^j_{\bull})$. That is,
$$\dim(V,W,\me)=\dim(\Gr(r,n))-\sum_{j=1}^s \codim(\omega_{I^j})$$
$$= r(n-r)-\sum_{j=1}^s\sum_{a=1}^r (n-r+a-i^j_a).$$
\end{defi}

\subsection{Induced flags}\label{inducedflags}
Suppose that $W$ is an $n$-dimensional vector space and $V\subseteq
W$ an $r$-dimensional subspace. Let $E_{\bull}$ be a complete flag
on $W$. This induces a complete flag on $V$ and a complete flag on
$W/ V$ by intersecting $E_{\bull}$ with $V$ and by projection
$p:W\to W/ V$  respectively. We denote these by $E_{\bull}(V)$ and
$E_{\bull}(W/ V)$ respectively. Explicitly, if $V\in
\Omega^o_I(E_{\bull})$ and $[n]\smallsetminus
I=\{\alpha(1)<\dots<\alpha(n-r)\}$, then $E_a(V)=E_{i_a}\cap V,\
a=1,\dots,r$ and $E_b(W/ V)=p(E_{\alpha(b)})\text{, }b=1,\dots,n-r.$
Given an ordered collection of flags $\me\in \Fl(W)^s$ we obtain
ordered collections of flags $\me(V)\in \Fl(V)^s$ and $\me(W/ V)\in
\Fl(W/ V)^s$ by performing the above operations in each coordinate
factor.

The following lemma follows from a direct calculation (see
~\cite{ful2}, Lemma 2 (i)).
\begin{lemma}\label{falcon}
Let $W$ be an $n$-dimensional vector space. Suppose $E_{\bull}\in
\Fl(W)$ and  $S\subseteq   V\subseteq   W $ are subspaces with
$\rk(V)=r$ and $\rk(S)=d$. Let $I$ be the unique subset of $[n]$ of
cardinality $r$  such that $V\in \Omega^o_{I}(E_{\bull})\subseteq
\Gr(r,W)$, and $K$ the unique subset of $[r]$ of cardinality $d$
such that $S\in \Omega^o_{K}(E_{\bull}(V))\subseteq  \Gr(d,V)$. Set
$L=\{i_a\mid a\in K\}$. Then,
 $S\in\Omega^o_L(E_{\bull})\subseteq  \Gr(d,W).$
\end{lemma}
\begin{lemma}\label{rappel}
Let $V$ be a $r$-dimensional subspace of an $n$-dimensional vector
space $W$, and $\me\in \Fl(W)^s$. Let $I^1,\dots, I^s$ be the unique
subsets of $[n]$ each of cardinality $r$ such that $V\in
\cap_{j=1}^s\Omega^o_{I^j}(E^j_{\bull})$. Suppose $S\subseteq V$ is
a $d$-dimensional vector subspace. Let $K^1,\dots, K^s$ be the
unique subsets of $[r]$ each of cardinality $d$ such that $S\in
\cap_{j=1}^s\Omega^o_{K^j}(E^j_{\bull}(V))$.

Then
$$\dim(S,V,\me(V))-\dim(S,W,\me)= \sum_{j=1}^s\sum_{a\in K^j} (n-r+a-i^j_a)-d(n-r).$$
\end{lemma}

\subsection{Tangent spaces}\label{tangentspaces}
 Let $V\in \Gr(r,W)$. If $\me\in
\Fl(W)^s$ and $I^1,\dots,I^s$ are the unique subsets of $[n]$ each
of cardinality $r$ such that $V\in \cap_{j=1}^s
\Omega^o_{I^j}(E^j_{\bull})$,  then the  tangent space at $V$ to the
scheme theoretic intersection
$\cap_{j=1}^s\Omega^o_{I^j}(E^j_{\bull})$ is given by (see
~\cite{gh})
\begin{equation}\label{tangent}
\{\phi \in \home(V,W/V)\mid \phi(E^j_a(V)) \subseteq   E^j_{i^j_a
-a}(W/V) \text{ for }a=1,\dots,r,\text{ }j=1,\dots,s\}.
\end{equation}
\begin{definition}\label{deftwo}
Let $\mi=(I^1,\dots,I^s)$ be a $s$-tuple of subsets of $[n]$ of
cardinality $r$ each. Let $V$  and $Q$ be vector spaces of rank $r$
and $n-r$ respectively and $(\mf,\mg)\in \Fl(V)^s\times \Fl(Q)^s$.
Define
\begin{equation}\label{intersection}
\home_{\mi}(V,Q,\mf,\mg)=\bigcap_{j=1}^s \{\phi \in \home(V,Q)\mid
\phi(F^j_a) \subseteq   G^j_{i^j_a -a}\text{ for }a=1,\dots,r\}.
\end{equation}
\begin{equation}\label{vs}
=\{\phi \in \home(V,Q)\mid \phi(F^j_a) \subseteq   G^j_{i^j_a
-a}\text{ for }a=1,\dots,r,\text{ }j=1,\dots,s\}
\end{equation}
\end{definition}
\subsection{Stratification and Universal families}\label{hamlett}
Let $\mi=(I^1,\dots,I^s)$, $V$ and $Q$ be as before. Let
$\mk=(K^1,\dots,K^s)$ be a $s$-tuple of subsets of $[r]$ each of
cardinality $d$. We consider the following ``universal objects''
(see ~\cite{gh} for more details).
\begin{enumerate}
\item[(A)] Define  $\mathfrak{H}_{\mi}(V,Q,\mk)$  to be the scheme
over  $\Fl(V)^s\times \Fl(Q)^s$ whose fiber over $(\mf,\mg)$ is
\begin{equation}\label{goo}
\{\phi\in \hoo\mid \rk(\ker(\phi))=d,\text{ }\ker(\phi)\in
\cap_{j=1}^s\Omega^o_{K^j}(F^j_{\bull})\}.
\end{equation}
\item[(B)] Define $\muu_{\mk}(V)$ to be the scheme over $\Fl(V)^s$ whose fiber over $\mf\in \Fl(V)^s$ is  $\cap_{j=1}^s\Omega^o_{K^j}(F^j_{\bull})\subseteq  \Gr(d,V)$.
\end{enumerate}
\begin{proposition}\label{labour}
\begin{enumerate}
\item $\mhh_{\mi}(V,Q,\mathcal{K})$  and $\mukv$ are smooth  and irreducible schemes.
\item If $\ham\neq\emptyset$, the natural morphism $p:\ham\to \mukv$ which maps $(\phi,\mf,\mg)$ to $(\ker(\phi),\mf)$, is smooth and surjective.
\item The dimension of $\mukv$ is $\dim \Fl(V)^s\ + \dim\Gr(d,r)\ -\ \sum_{j=1}^s \codim(\omega_{K^j}).$
\item The dimension of $\ham$ is
$$\dim\mukv + \dim \Fl(Q)^s + \{\sum_{j=1}^s\sum_{a\in K^j}(n-r+a-i^j_a)-d(n-r)\}.$$
\end{enumerate}
\end{proposition}

\section{The setting for  the proof of Fulton's
conjecture}\label{shostakovich}
 A subset $I\subseteq [n]$ of
cardinality $r$ also defines an irreducible representation
$V_{\lambda(I)}$ of $\GL(r)$, where
$\lambda(I)=(\lambda_1\geq\dots\geq \lambda_r)$ with
$\lambda_a=n-r+a-i_a$ for $a=1,\dots r$.

Now, let $I^1,\dots,I^s$ be subsets of $[n]=\{1,\dots,n\}$ each of
cardinality $r$. Assume that
$$\sum_{j=1}^{s}\sum_{a=1}^{r} (n-r +a -i^j_a)= r(n-r).$$
Let $\lambda_j=\lambda(I^j)$ be the weights of the corresponding
irreducible representations of $\GL(r)$.

Then it is well known that,
$$\prod_{j=1}^s\omega_{I^j}= \dim_{\Bbb{C}} [V_{\lambda_1}\tensor\dots\tensor V_{\lambda_s}]^{\SL(r)}[\text{ class of a point }]\in H^{*}(\Gr(r,n)).$$

The above notation and assumptions  will be kept fixed throughout
this paper and will be called the Fixed Setting. There is a related
space of parabolic vector spaces relevant to this setting:

\subsection{Parabolic vector spaces}
 A parabolic vector space $\tilde{V}$ is a $3$-tuple $(V,\mf,w)$, where $V$ is a vector space, $\mf\in\Fl(V)^s$ and $w$ is a function
$$w:\{1,\dots,s\}\times \{1,\dots,\rk(V)\}\to \Bbb{Z}$$ such that if
we let $w^j_l=w(j,l)$, the following holds for each $j=1,\dots,s$:
$$w^j_1\geq w^j_2\geq\dots\geq w^j_{\rk(V)}.$$

An isomorphism between parabolic vector spaces of the same rank and
weights $w$, $(V,\mf,w)$ and $(T,\mg,w)$ is an isomorphism $V\to T$
such that for any $j\in 1,\dots,s$ and $a<\rk(V)$ such that
$w^j_a>w^j_{a+1}$, $\phi(F^j_a)=G^j_{a}$. So in reality one ignores
the parts of the flags where the weights do not jump (we may
similarly define morphisms between parabolic vector spaces, and
create a corresponding abelian category).

Let $S\subseteq   V$ be a non zero subspace of rank $d$.  Let
$K^1,\dots,K^s$ be the unique subsets of $[\rk(V)]$ each of
cardinality $d$ such that $S\in
\cap_{j=1}^s\Omega^o_{K^j}(F^j_{\bull})$. Define the  parabolic
slope
$$\mu(S,\tilde{V})=\frac{\sum_{j=1}^s\sum_{a\in K^j} w^j_a}{d}.$$

A parabolic vector space $\tilde{V}$ is said to be semistable if for
each subspace $S\subseteq   V$, $\mu(S,\tilde{V})\leq
\mu(V,\tilde{V}).$

Given the fixed setting we get a choice of weights for parabolic
vector spaces. Here we consider parabolic vector spaces of the form
$(V,\mf,w)$ with $\rk(V)=r$ and
$$w^j_a=\lambda^j_a=n-r+a-i^j_a,\ j=1,\dots,s,\ a=1, \dots, r.$$
\subsection{Moduli spaces}
Let
$$\mm=\mathcal{M}(\mi,r,n)=\operatorname{Proj}\bigoplus_{N=1}^{\infty}(V^*_{N\lambda_1}\tensor\dots\tensor V^*_{N\lambda_s})^{\SL(r)}.$$
 be the (projective and irreducible)  moduli space of semistable parabolic vector spaces  with the
above weights. The proof of the properties below follows similar
properties for parabolic vector bundles.

Let $V$ be a $r$ dimensional vector space. There is an open subset
$U$ of $FL(V)^s$ formed by points $\mf$ so that $(V,\mf,w)$ is
semistable. There is a natural surjective map $U\to \mm$, and there
is a natural line bundle $\ml$ on $\mm$ obtained by descent via
Kempf's theory (see ~\cite{dmj}) of the natural line bundle
$\tilde{\ml}=\ml_{\lambda_1}\tensor\dots\tensor\ml_{\lambda_s}$ on
$\Fl(V)^s$ whose global sections are
$V^*_{\lambda_1}\tensor\dots\tensor V^*_{\lambda_s}$. In fact, note
that global sections (for example ~\cite{teleman})
\begin{equation}\label{dmitri}
H^0(\mm,\ml)=H^0(FL(V)^s,\tilde{\ml})^{\SL(V)}.
\end{equation}

 Similarly,
  \begin{equation} H^0(\mm,\ml^N)\ =\
[V^*_{N\lambda_1}\tensor\dots\tensor V^*_{N\lambda_s}]^{\SL(r)}.
\end{equation}
\subsection{Formulation of theorems as properties of $\mm$}
The saturation theorem of Knutson and Tao can therefore be
formulated as: For any positive integer $N$,
$$h^0(\mm,\ml)\neq 0 \Leftrightarrow h^0(\mm,\ml^N)\neq 0$$

Fulton's conjecture (theorem of Knutson, Tao and Woodward) can be
formulated as : For any positive integer $N$,
$$h^0(\mm,\ml)=1 \Leftrightarrow h^0(\mm,\ml^N)=1$$

In view of the ampleness of $\mathcal{L}$ on $\mm$, and the
connectedness of $\mm$ we have a reformulation of saturation and
Fulton's conjecture, the saturation statement is
$$\mm\neq \emptyset\Rightarrow h^0(\mm,\ml)\neq 0.$$

Fulton's conjecture is then the statement
$$\mm\neq \emptyset,\ \dim(\mm)>0 \Rightarrow h^0(\mm,\ml)\neq 1.$$

(This is applied to $\ml$ and $\ml^N$, the different linearization
$\tilde{\ml}^N$ does not change $\mm$ but changes the basic line
bundle on it to $\ml^N$.)
\subsection{The starting point}\label{SP}
Assume $\mm\neq \emptyset,\ \dim(\mm)>0$ and $h^0(\mm,\ml)=1$, and
 we will show that this leads to a contradiction. We will indicate the
starting point of the argument now. Let $\aA\in\ h^0(\mm,\ml)$ be
the unique non-vanishing section (upto scalars). Since $\ml$ is
ample and $\mm$ positive dimensional, {\bf the zero set of $\aA$ is
non-empty}. therefore $\aA$ vanishes at a point of the form
$(V,\mf)$ which is {\bf semistable}.

To ``test'' the assumption $h^0(\mm,\ml)=1$, we will need a  way of
producing global  sections of $\ml$. In ~\cite{invariant} we showed
a way of producing all sections of $h^0(\mm,\ml)$ via Schubert
calculus.
\subsection{Construction of sections of $\ml$}
We return to the notation of the fixed setting. Let $Q$ be a vector
space of rank $n-r$ and $\mg\in \Fl(Q)^s$. In ~\cite{invariant} we
showed that the pair $(Q,\mg)$ can be used to produce a section
$\aA(Q,\mg)$ of
$$H^0(\mm,\ml)\ =\ [V^*_{\lambda_1}\tensor\dots\tensor V^*_{\lambda_s}]^{\SL(r)}$$
We briefly recall the description: The zero set of $\aA(Q,\mg)$ on
$FL(V)^s$ is the set of points $(V,\me)$ for which the vector space
$\hoof$ (see Section ~\ref{tangentspaces}) is non-zero (this
condition is converted into a determinantal condition and hence a
section of the desired bundle.) If $W$ is an $n$-dimensional vector
space, $\mf\in\Fl(W)^s$ generic and
$\{V_1,\dots,V_m\}=\bigcap_{j=1}^s\Omega^o_{I^j}(F^j_{\bull})$, then
the sections $\aA(W/V_{\ell},\mf(W/V_{\ell}))$ give a basis for
$[V^*_{\lambda_1}\tensor\dots\tensor V^*_{\lambda_s}]^{\SL(r)}$ (in
fact $\aA(W/V_{\ell},\mf(W/V_{\ell}))$ vanishes at $(V_k,\mf(V_k))$
if and only if $\ell\neq k$)(ee Section ~\ref{inducedflags}).
\section{Return to the proof of Fulton's conjecture}
We now return to the situation at the end of Section ~\ref{SP}. Let
$Q$ be a vector space of rank $n-r$. Let $\aA\in
H^0(\mathcal{M},\ml)$ be the unique non-zero section (upto scalars).

{\bf Let $\mfz\subset FL(V)^s$ be the closure of an irreducible
component of the zero set of $\aA$ which contains a semistable
point} ($\mfz$ is to be fixed once and for all). Recall that the set
of semistable is open in any family (and by Equation ~\ref{dmitri},
we can consider sections of $\ml$ as invariant sections of
$\tilde{\ml}$ on $FL(V)^s$).

Let $(\me,\mh)$ be a generic point of $\mfz\times \Fl(Q)^s$. We know
that section $\aA(Q,\mh)$ is the unique non-zero section of $\ml$
upto scalars.  Therefore $\aA(Q,\mh)$ vanishes at $(V,\me)$, in
other words $\hoof$ is a non-zero vector space .  Let $\phi$ be a
generic element of $\home_{\mi}(V,Q,\me,\mh)$. Let $S=\ker(\phi)$,
$d=\rk(S)$ ($d=0$ is possible!), and let $\mk=(K^1,\dots, K^s)$ be
the unique $s$-tuple of subsets of $[r]$ each of cardinality $d$
such  that  $S\in \cap_{j=1}^s\Omega^o_{K^j}(E^j_{\bull})$.
\begin{proposition}\label{one}
Let $(V,\mf)$ be such that
$\cap_{j=1}^s\Omega^o_{K^j}(F^j_{\bull})\neq \emptyset$. Then
$\mf\in \mfz$.
\end{proposition}
\begin{proof}
We consider the spaces $\ham$ and $\mukv$ from ~\cite{gh} (and
recalled in Section ~\ref{hamlett} of this paper):

Let $X$ be the nonempty open subset of $\Fl(Q)^s$ formed by points
$\mg$ so that $\theta(Q,\mg)\neq 0$. We claim that if
$$Y=\{(V,\mf,\mg)\in\ham\mid \mg\in X\}$$
 and $(V,\mf,\mg)\in Y$ then $\aA$ vanishes at $(V,\mf)$. This is clear because $\aA(Q,\mg)$ vanishes at $(V,\mf)$ and $\aA$ is a multiple of $\aA(Q,\mg)$.

Clearly $Y$ is non-empty (by assumption!) and dense in the
irreducible $\ham$. Therefore $\aA$ vanishes at all points of the
form $(V,\mf)$ for which there is a  point of the form
$(V,\mf,\mg)\in\ham$. The surjectivity of $\ham\to \mukv$ therefore
assures us that $\aA$ vanishes at $(V,\mf)$ if there exists a point
of the form $(S,\mf)\in \mukv$. Since $\mukv$ is irreducible, the
proof is complete.
\end{proof}
Let $\mfa\subseteq\Fl(V)^s$ be the closure of the image of  $\mukv$.
By Proposition ~\ref{one} $\mfa\subseteq \mfz$ (note that $\mfa$
being the closure of the image of $\mukv$ is irreducible). By
assumption $\mfz\subseteq \mfa$. Hence
\begin{lemma}
$\mfa=\mfz$.
\end{lemma}

Let $R,T$ be vector spaces of dimension $d$ and $r-d$ respectively,
and $U$ a nonempty open subset of  $\Fl(R)^s\times \Fl(T)^s$ which
is stable under $\GL(R)^s\times\GL(T)^s$ (so that one has
canonically defined open subsets of $\Fl(R')^s\times\Fl(T')^s$ for
any vector spaces $R'$ and $T'$ of ranks $d$ and $r-d$
respectively).

\begin{proposition}\label{fulton}
 There exists a nonempty open subset $\tilde{U}$ of $\mfz\times\Fl(Q)^s$ such that
 for $(\mf,\mg)\in \tilde{U}$,
\begin{enumerate}
\item[(a)] The intersection $\cap_{j=1}^s\Omega^o_{K^j}(F^j_{\bull})$  is equidimensional of dimension $\dim \mukv -\dim \mfz$.
\item[(b)] If $\phi$ is a general element of $\hoo$, and $S=\ker(\phi)$ then the induced pair of flags $(\mf(S),\mf(V/S))$ ``is a point of $U$''.
\item[(c)] The rank of $\hoo$ is
$$\dim\mukv -\dim \mfz +\  \{\sum_{j=1}^s\sum_{a\in K^j}(n-r+a-i^j_a)-d(n-r)\}.$$
\end{enumerate}
\end{proposition}
\begin{proof}\
Item (a) follows from generic flatness of $\mukv\to \mfa=\mfz$.

Let $W_1\subseteq \mukv$ be the non-empty open subset of points
$(S,\mf)$ such that  $(\mf(S),\mf(V/S))$ is a point of $U$. Let
$W_2$ be the inverse image of $W_1$ in $\ham$. Let $\tilde{U}$ be an
open subset of $\mfz\times\Fl(Q)^s$ so that the map $[\ham-W_2]\to
\mfz\times\Fl(Q)^s$ is flat over $\tilde{U}$. This proves (b).

The fiber dimension of $\ham\to \mfz\times\Fl(Q)^s$ is easily seen
to be
$$\dim \mukv-\dim \mfz + \{\sum_{j=1}^s\sum_{a\in K^j}(n-r+a-i^j_a)-d(n-r)\}.$$
This proves (c).
\end{proof}
\subsubsection{The conclusion of the  proof of Fulton's conjecture }
Let $(\mf,\mg)$ be a general point of $\mfz\times \Fl(Q)^s$ as in
Proposition ~\ref{fulton}. Let $\phi$ be a general point of $\hoo$,
$S=S(1)=\ker(\phi)$. Then the induced flags $(S,\mf(S))$ and
$(V/S,\mf(V/S))$ can be assumed to be general and in mutually
general position This follows from  proposition ~\ref{fulton}.

Set $L^j(1)=K^j$ for $j=1,\dots,s$. Now suppose that
$\cap_{j=1}^s\Omega^o_{K^j}(F^j_{\bull})$ is positive dimensional at
$S=S^{(1)}$.

If $\psi^{(1)}$ is a generic element in the tangent space of
$\cap_{j=1}^s\Omega^o_{K^j}(F^j)\subseteq   \Gr(d,V)$ at $S$, we can
view $\psi^{(1)}$ as a map $S\to V/S$ (the tangent space to
$\Gr(d,V)$ at $S$ is $\home(S,V/S)$). Let $S^{(2)}$ be the kernel of
$\psi^{(1)}$ and assume that
 $S^{(2)}\in \cap_{j=1}^s\Omega^o_{L^j(2)}(F^j(S))\subseteq   \Gr(d,S)$.

We  obtain a sequence of inclusions
$$S^{(h)}\subsetneq S^{(h-1)}\subsetneq\dots\subsetneq S^{(1)}=S\subsetneq V$$ inductively as follows. Assume
 $$S^{(\ell)}\in \cap_{j=1}^s\Omega^o_{L^{j}(\ell)}(F^j(S^{(\ell-1)}))\subseteq   \Gr(d,S^{(\ell-1)}).$$ If this intersection is $0$ dimensional at $S^{(\ell)}$then the process stops at $h=\ell$. If it positive dimensional at $S$,
 let $\psi^{(\ell)}$ be the generic element of the tangent space at $S^{(\ell)}$ of
$$\cap_{j=1}^s\Omega^o_{L^j(\ell)}(F^j(S^{(\ell-1)}))\subseteq  \Gr(d,S^{(\ell-1)}).$$
View $\psi^{(\ell)}$ as a map $S^{(\ell)}\to
S^{(\ell-1)}/S^{(\ell)}$ and define $S^{(\ell+1)}=\ker
\psi^{(\ell)}$. And continue with the ``recursion''. This procedure
starting from $S$ will be called the ``tangent space method''.

For $u=1,\dots,h$, let $d_u$ be the rank of $S^{(u)}$,
${\Ftwo}(u)=({\Fone}^1(u),\dots,{\Fone}^s(u))$ the unique $s$-tuple
of subsets of $[r]$ each of cardinality $d_u$ such that  $S^{(u)}\in
\cap_{j=1}^s\Omega^o_{{\Fone}^j(u)}(F^j_{\bull})\subseteq
\Gr(d_u,V)$,
 By Lemma ~\ref{falcon} (applied to $S^{(u)}\subseteq S\subseteq V$ and $\mf\in \Fl(V)^s$)
\begin{equation}\label{r1emember}
{\Fone}^j(u)=\{k^j_b\mid b\in L^j(u)\}
\end{equation}
We claim
\begin{proposition}\label{fin}
\begin{enumerate}
\item[(i)] $\dim\mukv-\dim \mfz=\dim \cap_{j=1}^s\Omega^o_{K^j}(F^j_{\bull})$  is less than or equal to (in fact equal to, we will not need this)
$$\dim(S,V,\mf) +\dim(S^{(h)},S,\mf(S))-\dim (S^{(h)},V,\mf).$$
\item[(ii)]
$$\dim(S,V,\mf)+\{\sum_{j=1}^s\sum_{a\in K^j}(n-r+a-i^j_a)-d_u(n-r)\}$$
$$\leq
\dim(S^{(h)},V,\mf)- \dim(S^{(h)},S,\mf(S))
+\{\sum_{j=1}^s\sum_{a\in {\Fone}^j(h)}(n-r+a-i^j_a)-d_h(n-r)\}
$$
\end{enumerate}
\end{proposition}
From (i) and (ii), and Proposition ~\ref{fulton},  we conclude that
the rank of $\hoo$ is no greater than
\begin{equation}\label{clincher}
\sum_{j=1}^s\sum_{a\in {\Fone}^j(h)}(n-r+a-i^j_a)-d_h(n-r)
\end{equation}

The semistability of $(V,\mf,w)$ we obtain that Expression
~\ref{clincher}
 is $\leq 0$ and hence the geometric proof of Fulton's conjecture would be complete once Proposition ~\ref{fin} is proved.

\begin{proof}(Of Proposition ~\ref{fin})
The dimension of $\dim \cap_{j=1}^s\Omega^o_{K^j}(F^j_{\bull})$ is
no more than the dimension of its tangent space at $S$ which is
$\home_{\mk}(S,V/S,\mf(S),\mf(V/S))$ and using the main theorem in
~\cite{gh} (and Lemma ~\ref{rappel}) we find that (i) holds.

Item (ii) follows from the filtration lemma in ~\cite{gh} (recalled
below) where we apply it with $\eta_u=\phi\circ \psi^{(u)}$.
\end{proof}
\subsection{The Filtration Lemma}
For a vector space  $W$ of rank $n$,  define  $B(W)\subseteq
\Fl(W)^s$ to be the largest Zariski open subset of $\Fl(W)^s$
satisfying the following property: If $\me\in B(W)$ and
$\mi=(I^1,\dots,I^s)$ a $s$-tuple of subsets of $[n]$ each of the
same cardinality $r$, then every irreducible component of the
intersection $\cap_{j=1}^s\Omega_{I^j}(E^j_{\bull})$ (which is
possibly empty) is proper. By Kleiman's transversality theorem
~\cite{kl}, it follows that $B(W)$ is nonempty.

\begin{lemma} Consider a $5$-tuple of the form $(V,Q,\mf,\mg,\mi)$ where $V$ and $Q$ are non-zero vector spaces of ranks $r$ and $n-r$ respectively, $\mi=(I^1,\dots,I^s)$ a $s$-tuple of subsets of $[n]$ each of cardinality $r$ and   $\mg\in B(Q)$.

Suppose in addition that we are given a filtration by vector
subspaces
\begin{equation}\label{filtratione}
S^{(h)}\subsetneq S^{(h-1)}\subsetneq\dots \subsetneq
S^{(1)}\subsetneq S^{(0)}=V
\end{equation}
and injections (of vector spaces) from the graded quotients
$\eta_{u}:{S^{(u)}}/{S^{(u+1)}}\hookrightarrow Q$ for $u=0,\dots,
h-1$ such that for  $j=1,\dots,s$ and $a=1,\dots,r$, (where we write
$\eta_u$ again  for the composite $S^{(u)}\to
{S^{(u)}}/{S^{(u+1)}}\leto{\eta_u} Q$)
$$\eta_u(S^{(u)}\cap F^j_a)\subseteq  G^j_{i^j_a-a}$$

Then, letting $d_u$ be the rank of $S^{(u)}$,
${\Ftwo}(u)=({\Fone}^1(u),\dots,{\Fone}^s(u))$ the unique $s$-tuple
of subsets of $[r]$ each of cardinality $d_u$ such that  $S^{(u)}\in
\cap_{j=1}^s\Omega^o_{{\Fone}^j(u)}(F^j_{\bull})\subseteq\Gr(d_u,V)$
and
$$\dim(S,V,\mf)+\{\sum_{j=1}^s\sum_{a\in K^j}(n-r+a-i^j_a)-d_u(n-r)\}$$
$$\leq
\dim(S^{(h)},V,\mf)- \dim(S^{(h)},S,\mf(S))
+\{\sum_{j=1}^s\sum_{a\in {\Fone}^j(h)}(n-r+a-i^j_a)-d_h(n-r)\}
$$
\end{lemma}
\appendix\section{Resume of results in ~\cite{gh}}
Let $\mi=(I^1,\dots,I^s)$ be a $s$-tuple of subsets of $[n]$ of
cardinality $r$ each and $W$ an n-dimensional vector space. Let us
ask the following questions:
\begin{enumerate}
\item[(Q1)] For generic $\me\in\Fl(W)^s$  is $\cap_{j=1}^s\Omega^o_{I^j}(E^j_{\bull})$ empty?
\item[(Q2)] Let $\mfa\subseteq  \Fl(V)^s$ be the set of $\me$ such that $\cap_{j=1}^s\Omega^o_{I^j}(E^j_{\bull})$ is non-empty. For generic ($\mfa$ is the image of $\muiw$ and is hence  irreducible) $\me\in \mfa$ what is then the dimension of (each irreducible component of) $\cap_{j=1}^s\Omega^o_{I^j}(E^j_{\bull})$?
\end{enumerate}

 Let $V$  and $Q$ be vector spaces of rank $r$ and $n-r$ respectively and $(\mf,\mg)\in \Fl(V)^s\times \Fl(Q)^s$ a generic point. Then, to answer these questions, according to ~\cite{gh} one needs to proceed as follows. The answer for (Q2) is the rank of $\hoo$. If this rank equals the expected dimension $[\dim(Gr(r,n))-\sum_{j=1}^s \codim(\omega_{I^j})]$, then the answer to (Q1) is affirmative (and vice-versa).
\begin{theorem}
There exists a   filtration by vector subspaces obtained by the
``tangent space'' method
\begin{equation}\label{filtratione}
S^{(h)}\subsetneq S^{(h-1)}\subsetneq\dots \subsetneq
S^{(1)}\subsetneq S^{(0)}=V
\end{equation}
and injections (of vector spaces) from the graded quotients
$\eta_{u}:{S^{(u)}}/{S^{(u+1)}}\hookrightarrow Q$ for $u=0,\dots,
h-1$ such that the following property is satisfied: For
$u=1,\dots,h$, let $d_u$ be the rank of $S^{(u)}$,
${\Ftwo}(u)=({\Fone}^1(u),\dots,{\Fone}^s(u))$ the unique $s$-tuple
of subsets of $[r]$ each of cardinality $d_u$ such that  $S^{(u)}\in
\cap_{j=1}^s\Omega^o_{{\Fone}^j(u)}(F^j_{\bull})\subseteq
\Gr(d_u,V)$, then
\begin{enumerate}
\item[(i)] $\dim(S^{(h)},V,\mf)=0.$
\item[(ii)] For $u=0,\dots,h-1$, $j=1,\dots,s$ and $a=1,\dots,r$, (where we write $\eta_u$ again  for the composite $S^{(u)}\to {S^{(u)}}/{S^{(u+1)}}\leto{\eta_u} Q$)
$$\eta_u(S^{(u)}\cap F^j_a)\subseteq  G^j_{i^j_a-a}$$
\item[(iii)] The vector space $\hoo$ is of rank (the second term of the expression below is the same as the quantity appearing in Inequality $(\dagger_{{\Ftwo}(h)}^{\mi})$)
\begin{equation}\label{une}
[\dim(Gr(r,n))-\sum_{j=1}^s \codim(\omega_{I^j})]+
\{\sum_{j=1}^s\sum_{a\in {\Fone}^j(h)}(n-r+a-i^j_a)-d_h(n-r)\}.
\end{equation}
\end{enumerate}
\end{theorem}
\begin{remark}
The filtration is constructed in the course of the proof of this
theorem. Here we start with a generic element of $\hoo$, let
$S=\ker(\phi)\subset V$ and apply the tangent space method to it.
\end{remark}

\bibliographystyle{plain}
\def\noopsort#1{}

\end{document}